\documentclass[reqno]{udt14a_sciendo_CC2} 

\usepackage{hyperref} 
\makeatletter \thm@headfont{\bfseries\scshape} \makeatother
\usepackage{latexsym,euscript}
\usepackage{amsbsy,amscd,amsfonts,amsmath,amsopn,amssymb,amstext,amsthm,amsxtra,mathtools,bbm,mathrsfs}
\usepackage{epsf,epsfig,graphics,color,verbatim}
\usepackage{cite, float,multirow}
\usepackage{subfigure}
\usepackage{enumerate}
\theoremstyle{plain}
\newtheorem{theorem}{Theorem}[section]
\newtheorem{lemma}[theorem]{Lemma}
\newtheorem{corollary}[theorem]{Corollary}

\theoremstyle{definition} 
\newtheorem{definition}[theorem]{Definition}

\newtheorem{remark}{Remark}

\newtheorem*{acknowledgement}{Acknowledgement}

\pagespan{1}{2}

\begin{document}

\title[A generalization of van der Corput's Difference Theorem]
{A generalization of van der Corput's Difference Theorem}

\author[1,2]{Sohail Farhangi
        }
\affil[1]{Adam Mickiewicz University, Pozna\'n, Poland}
\address{\bf{Sohail Farhangi}\\
         Department of Mathematics and Informatics        \\
         University of Adam Mickiewicz      \\
         ul. Wieniawskiego 1        \\
         61-712--Pozna\'n \\
         Poland}
\email{sohail.farhangi@gmail.com}
\def\shortauthors{S. Farhangi}               

\keywords{van der Corput's Difference Theorem, discrepancy, ergodic hierarchy of mixing properties, }
\thanks {Supported by the Grant “Set theoretic methods in dynamics and number theory”, n. 2019/34/E/ST1/00082.}
\begin{abstract}
We prove a generalization of van der Corput's Difference Theorem in the theory of uniform distribution by establishing a connection with unitary operators that have Lebesgue spectrum. This allows us to show, for example, that if $(x_n)_{n = 1}^\infty \subseteq [0,1]$ is such that $(x_{n+h}-x_n)_{n = 1}^\infty$ is uniformly distributed for all $h \in \mathbb{N}$, then $(x_{n_k})_{k = 1}^\infty$ is uniformly distributed, where $(n_k)_{k = 1}^\infty$ is an enumeration of the $1s$ in the classical Thue-Morse sequence. We also establish a variant of van der Corput's Difference Theorem that is connected to unitary operators with continuous spectrum. Lastly, we obtain a new characterization of those sequence $(x_n)_{n = 1}^\infty \subseteq [0,1]$ for which $(x_{n+h},x_n)_{n = 1}^\infty$ is uniformly distributed in $[0,1]^2$ for all $h \in \mathbb{N}$. 
\end{abstract}

\maketitle

\par 

\section{Introduction}
In \cite{OriginalvanderCorput} van der Corput proved Theorem \ref{UniformDistributionvanderCorput'sDifferenceTheorem}, which is now known as van der Corput's Difference Theorem (henceforth abbreviated as vdCDT).

\begin{theorem}[{\cite[Theorem 1.3.1]{Kuipers&Niederreiter}}]
\label{UniformDistributionvanderCorput'sDifferenceTheorem}
If $(x_n)_{n = 1}^{\infty} \subseteq [0,1]$ is a sequence for which $(x_{n+h}-x_n)_{n = 1}^{\infty}$ is uniformly distributed for all $h \in \mathbb{N}$, then $(x_n)_{n = 1}^{\infty}$ is uniformly distributed.
\end{theorem}

It is clear that there exist uniformly distributed sequences $(x_n)_{n = 1}^\infty \subseteq [0,1]$ for which $(x_{n+h}-x_n)_{n = 1}^\infty$ is not uniformly distributed for any $h \in \mathbb{N}$. A natural example of such a sequence is found by taking some $\alpha \in \mathbb{R}\setminus\mathbb{Q}$ and letting $x_n = I_n := \lfloor\sqrt{n}\rfloor\alpha$, and it is worth noting that in this example, for any $h \in \mathbb{N}$ we have $x_{n+h}-x_n = h\alpha$ for a full density set of $n \in \mathbb{N}$. It is therefore natural to ask what additional properties are satisfied by those $(x_n)_{n = 1}^\infty \subseteq [0,1]$ that satisfy the hypothesis of vdCDT, and it is the goal of the current paper to answer this question.

For reasons that will become apparant later, let us call a sequence $(x_n)_{n = 1}^\infty \subseteq [0,1]$ a \textbf{sL-sequence} (spectrally Lebesgue sequence) if $(x_{n+h}-x_n)_{n = 1}^\infty$ is uniformly distributed for all $h \in \mathbb{N}$. It had already been observed in \cite{vdCSetsK&MF} that if $(x_n)_{n = 1}^\infty$ is a sL-sequence, then so is $(x_n+n\alpha)_{n = 1}^\infty$ for any $\alpha \in \mathbb{R}$, and this can be used to show that $(x_{n_k})_{k = 1}^\infty$ is uniformly distributed when $(n_k)_{k = 1}^\infty \subseteq \mathbb{N}$ is almost periodic (see \cite[Theorem 4.4]{vanderCorputTheoremSurvey}), or that $(x_n,n\alpha)_{n = 1}^\infty$ is uniformly distributed in $[0,1]^2$ when $\alpha \in \mathbb{R}\setminus\mathbb{Q}$. However, we see that these properties are also satisfied by the sequence $(I_n)_{n = 1}^\infty$ that we constructed above, so they cannot possibly characterize sL-sequences. One main goal of this paper is to better understand sL-sequences through similar corollaries. In particular, we will give a sufficient condition on those $(n_k)_{k = 1}^\infty \subseteq \mathbb{N}$ for which $(x_{n_k})_{k = 1}^\infty$ is uniformly distributed whenever $(x_n)_{n = 1}^\infty$ is a sL-sequence, and we will give a sufficient condition on those $(y_n)_{n = 1}^\infty \subseteq [0,1]$ for which $(x_n,y_n)_{n = 1}^\infty$ is uniformly distributed in $[0,1]^2$ whenever $(x_n)_{n = 1}^\infty$ is a sL-sequence.

In order to give some context for our main results, let us recall the following Hilbertian analogues of Theorem \ref{UniformDistributionvanderCorput'sDifferenceTheorem} that were introduced by Bergelson in \cite{WMPET} and are of great use in Ergodic Theory.

\begin{theorem}
\label{ClassicalvanderCorput'sDifferenceTheorems}
Let $\mathcal{H}$ be a Hilbert space and $(x_n)_{n = 1}^{\infty} \subseteq \mathcal{H}$ a bounded sequence of vectors. 

\begin{enumerate}
    \item[(i)] If for every $h \in \mathbb{N}$ we have

    \begin{equation}
        \lim_{N\rightarrow\infty}\frac{1}{N}\sum_{n = 1}^N\langle x_{n+h}, x_n\rangle = 0\text{, then }\lim_{N\rightarrow\infty}\left|\left|\frac{1}{N}\sum_{n = 1}^Nx_n\right|\right| = 0.
    \end{equation}
    
    \item[(ii)] If

    \begin{equation}
        \lim_{h\rightarrow\infty}\limsup_{N\rightarrow\infty}\left|\frac{1}{N}\sum_{n = 1}^N\langle x_{n+h}, x_n\rangle\right| = 0\text{, then }\lim_{N\rightarrow\infty}\left|\left|\frac{1}{N}\sum_{n = 1}^Nx_n\right|\right| = 0.
    \end{equation}
    
    \item[(iii)] If
    
    \begin{equation}
        \lim_{H\rightarrow\infty}\frac{1}{H}\sum_{h = 1}^H\limsup_{N\rightarrow\infty}\left|\frac{1}{N}\sum_{n = 1}^N\langle x_{n+h}, x_n\rangle\right| = 0\text{, then }\lim_{N\rightarrow\infty}\left|\left|\frac{1}{N}\sum_{n = 1}^Nx_n\right|\right| = 0.
    \end{equation} 
\end{enumerate}
\end{theorem}

It is natural to ask why there are 3 variations of vdCDT when working with Hilbert spaces, but only 1 vdCDT in the Theory of Uniform Distribution. This question is addressed in \cite[Chapter 2]{SohailsPhDThesis} by establishing connections between variations of vdCDT with various levels of the ergodic hierarchy of mixing properties of a unitary operator. In particular, \cite[Corollary 2.2.11]{SohailsPhDThesis} shows that \ref{ClassicalvanderCorput'sDifferenceTheorems}(iii) corresponds to weak mixing, \cite[Corollary 2.2.15]{SohailsPhDThesis} shows that \ref{ClassicalvanderCorput'sDifferenceTheorems}(ii) corresponds to strong mixing, and \cite[Corollary 2.2.17]{SohailsPhDThesis} and \cite{LebesgueSpectrumvanderCorput} show that \ref{ClassicalvanderCorput'sDifferenceTheorems}(i) corresponds to Lebesgue spectrum. 

Building upon these results, \cite[Corollary 2.4.18]{SohailsPhDThesis} is a variation of Theorem \ref{UniformDistributionvanderCorput'sDifferenceTheorem} which resembles Theorem \ref{ClassicalvanderCorput'sDifferenceTheorems}(iii) and is connected to a new notion in the Theory of Uniform Distribution, wm-sequences, that is related to weak mixing. Similarly, \cite[Corollary 2.4.22]{SohailsPhDThesis} is a variation of Theorem \ref{UniformDistributionvanderCorput'sDifferenceTheorem} which resembles Theorem \ref{ClassicalvanderCorput'sDifferenceTheorems}(ii) and is connected to another new notion in the Theory of Uniform Distribution, sm-sequences, that is related to strong mixing. Interestingly, a natural generalization of Theorem \ref{UniformDistributionvanderCorput'sDifferenceTheorem} was not obtained in \cite[Chapter 2.4]{SohailsPhDThesis}, so it is one of the (aforementioned) goals of the present paper to prove such a result. Another main goal of this paper is to give new characterizations of wm-sequences. In particular, we identify the classes of c-sequences\footnote{We remark the the sequence $(I_n)_{n = 1}^\infty$ is a uniformly distributed c-sequence.} and compact sequences, which are related to almost periodicity. Then we show that $(x_n)_{n = 1}^\infty \subseteq [0,1]$ is a wm-sequence if and only if $(x_n,y_n)_{n = 1}^\infty$ is uniformly distributed in $[0,1]^2$ for any uniformly distributed c-sequences $(y_n)_{n = 1}^\infty \subseteq [0,1]$. Similarly, we show that $(x_n)_{n = 1}^\infty$ is a wm-sequence if and only if $(x_{n_k})_{k = 1}^\infty$ is uniformly distributed for any compact sequence $(n_k)_{k = 1}^\infty \subseteq \mathbb{N}$.

Lastly, we define a class of mixing sequences $(x_n)_{n = 1}^\infty \subseteq [0,1]$ called o-sequences, and we show that they are the same as those sequences for which $(x_{n+h},x_n)_{n = 1}^\infty$ is uniformly distributed in $[0,1]^2$ for all $h \in \mathbb{N}$. While the definition of an o-sequence appears to be intimately related to Theorem \ref{UniformDistributionvanderCorput'sDifferenceTheorem}, we give two examples of sL-sequences that are not o-sequences.

We must also mention here that variations of Theorem \ref{UniformDistributionvanderCorput'sDifferenceTheorem} and \ref{ClassicalvanderCorput'sDifferenceTheorems} corresponding to ergodicity and mild mixing were already established in \cite[Chapter 2]{SohailsPhDThesis}, and connections between other forms of vdCDT in Hilbert spaces and the ergodic hierarchy of mixing were established in \cite{UnifiedVDC}.

The structure of the paper is as follows. In the next Section \ref{PreliminariesSection} we gather facts discrepancy, spectral theory, and mixing/rigid sequence that we will need later on. In Section \ref{MainResultsSection} we state and prove our main results. In Section \ref{SectionWithExamples} we list some examples that give more context to our main results

\section{Preliminaries}\label{PreliminariesSection}

\subsection{Discrepancy}
\begin{definition}
\label{DefinitionOfDiscrepancy}
Let $d \in \mathbb{N}$. The \textbf{discrepancy} of $(x_n)_{n = 1}^N \subseteq [0,1]^d$ is given by

\begin{equation}
D_N\left((x_n)_{n = 1}^N\right) := \underset{B \in \mathcal{R}}{\text{sup  }}\left|\frac{1}{N}|\{1 \le n \le N\ |\ x_n \in B\}|-m^d(B)\right|,
\end{equation}
where $\mathcal{R}$ denotes the collection of all rectangular prisms contained in $[0,1]^d$. For an infinite sequence $(x_n)_{n = 1}^{\infty} \subseteq [0,1]^d$, we let

\begin{equation}
    \overline{D}((x_n)_{n = 1}^{\infty}) := \limsup_{N\rightarrow\infty}D_N((x_n)_{n = 1}^N),
\end{equation}
and we let 

\begin{equation}
D((x_n)_{n = 1}^{\infty},(N_q)_{q = 1}^{\infty}) := \lim_{q\rightarrow\infty}D_{N_q}((x_n)_{n = 1}^{N_q}),
\end{equation}
provided that the limit exists.
\end{definition}

It is worth noting that a sequence $(x_n)_{n = 1}^\infty \subseteq [0,1]^d$ is uniformly distributed if and only if $\overline{D}((x_n)_{n = 1}^\infty) = 0$ (cf. Theorem 2.1.1 in \cite{Kuipers&Niederreiter}). We will also be needing the discrepancy of $(n\alpha)_{n = 1}^N$ for some $\alpha$, so we record the following result that is a consequence of Theorem 2.3.2 and Example 2.3.1 of \cite{Kuipers&Niederreiter}.

\begin{theorem}\label{DiscrepancyOfRotationByAnAlgebraicIrrational}
    Let $\alpha \in \mathbb{R}$ be an algebraic irrational. For any $\epsilon > 0$, there exists a constant $C = C(\alpha,\epsilon)$ such that for all $N \in \mathbb{N}$ we have

    \begin{equation}
        D_N\left((x_n)_{n = 1}^N\right) \le CN^{-1+\epsilon}.
    \end{equation}
\end{theorem}

The following is a corollary of the Koksma-Hlawka inequality, and is a special case of \cite[Theorem 2.5.6]{Kuipers&Niederreiter}.

\begin{theorem}\label{SpecialCaseOfKHInequality}
    Let $d \in \mathbb{N}$, let $\vec{m} := (m_1,\cdots,m_d) \in \mathbb{Z}^d\setminus\{0\}$, and let $M = \max_{1 \le i \le d}(|m_i|)$. For all $(x_n)_{n = 1}^N \in [0,1]^d$ we have

    \begin{equation}
        \left|\frac{1}{N}\sum_{n = 1}^Ne(\vec{m}\cdot x_n)\right| \le (4\pi M)^dD_N\left((x_n)_{n = 1}^N\right)
    \end{equation}
\end{theorem}
\subsection{Spectral theory}\label{SpectralTheorySubsection}
We give a brief review of spectral theory in order to give the reader context for the definitions and results in later sections. For a more detailed treatment we refer the reader to \cite[Chapter 18]{OTAoET}.

Suppose that $\mathcal{H}$ is a Hilbert space, $U:\mathcal{H}\rightarrow\mathcal{H}$ is a unitary operator, and $\xi \in \mathcal{H}$ is a cyclic vector for $U$. The Spectral Theorem tells us that there is a Hilbert space isomorphism $\phi:H\rightarrow L^2([0,1],\nu_\xi)$ for which $\phi U = M_e\phi$, where $M_e:L^2([0,1],\nu_\xi)\rightarrow L^2([0,1],\nu_\xi)$ is the unitary operator given by $(M_ef)(x) = e(-x)f(x)$, and $\phi(\xi) = 1$. The measure $\nu_\xi$ is the spectral measure of $\xi$ (with respect to $U$) and is seen to satisfy

\begin{equation}
    \hat{\nu}_\xi(n) = \int_0^1e(-nx)d\nu_\xi(x) = \langle M_e^n1,1\rangle_{L^2} = \langle U^n\xi, \xi\rangle_\mathcal{H}.
\end{equation}
Furthermore, $\xi_1,\xi_2 \in \mathcal{H}$ are such that $\nu_{\xi_1}\perp\nu_{\xi_2}$, then $\langle \xi_1,\xi_2\rangle_\mathcal{H} = 0$. Some classical results in ergodic theory can be understood as a corollary of this fact. 

For example, if $\xi_1,\xi_2 \in \mathcal{H}$ are such that $U\xi_1 = \xi_1$ and $\xi_2 \in c\ell(\{U\eta-\eta\ |\ \eta \in \mathcal{H}\})$, then $\nu_{\xi_1} = ||\xi_1||\mathbbm{1}_{\{0\}}$ and $\nu_{\xi_2}(\{0\}) = 0$, hence $\nu_{\xi_1}\perp\nu_{\xi_2}$, and we recover the the decomposition of Von Neumann of $\mathcal{H}$ into the invariant space and the span closure of cocycles. Another example is when $\xi_1$ is in the span closure of eigenvectors of $U$ and $\xi_2$ is a weakly mixing vector for $U$, in which case $\nu_{\xi_1}$ is a discrete measure and $\nu_{\xi_2}$ is a continuous measure. We again see that $\nu_{\xi_1}\perp\nu_{\xi_2}$, and we recover the Jacobs-de Leeuw-Glicksberg decomposition for a unitary operator. Another example that we will be using later on is when $\nu_{\xi_1}\perp m$ and $\nu_{\xi_2} << m$, in which case we immediately see that $\nu_{\xi_1}\perp\nu_{\xi_2}$, hence $\langle \xi_1,\xi_2\rangle_\mathcal{H} = 0$.
\subsection{Mixing and rigidity of sequences}
\subsubsection{Definitions}
We begin with the notion of permissible triples and permissible pairs, which are notions used to study bounded sequences by considering all possible subsequences in which relevant limits exists. Permissible triples were used in \cite[Chapter 2]{SohailsPhDThesis} and \cite{LebesgueSpectrumvanderCorput} to construct Hilbert spaces associated to a given pair of sequences, and we record the results of this construction as Lemma \ref{WhyWeNeedPermissibleTriples}. The usage of permissible triples is similar to the construction of Furstenberg Systems discussed in \cite{FSystemsOfBddedMultFunctions} and the references therein. Lemma \ref{ExampleNecessitatingPermissibleTriples} (\cite[Lemma 2.2.20]{SohailsPhDThesis}) will indicate that it is indeed necessary to consider permissible pairs (triples) instead of some simpler notion.

\begin{definition}
    Let $(x_n)_{n = 1}^\infty$ and $(y_n)_{n = 1}^\infty$ be bounded sequences of complex numbers and let $(N_q)_{q = 1}^\infty \subseteq \mathbb{N}$ be an increasing sequence. The triple $((x_n)_{n = 1}^\infty,$ $(y_n)_{n = 1}^\infty,(N_q)_{q = 1}^\infty)$ is a \textbf{permissible triple} if

    \begin{equation}
        \lim_{q\rightarrow\infty}\frac{1}{N_q}\sum_{n = 1}^{N_q}c_{n+h}\overline{d_n}
    \end{equation}
    exists for all $h \in \mathbb{N}$ and $(c_n)_{n = 1}^\infty,(d_n)_{n = 1}^\infty \in \{(x_n)_{n = 1}^\infty,$ $(y_n)_{n = 1}^\infty\}$. Given a family $\{(x_{n,h})_{n = 1}^{\infty}\}_{h = 1}^{\infty}$ of sequences in $[0,1]^d$ and an increasing sequence $(N_q)_{q = 1}^{\infty} \subseteq \mathbb{N}$, we define $(\{(x_{n,h})_{n = 1}^{\infty}\}_{h = 1}^{\infty},(N_q)_{q = 1}^{\infty})$ to be a \textbf{permissible pair} if for all $h \in \mathbb{N}$ $D((x_{n,h})_{n = 1}^{\infty},(N_q)_{q = 1}^{\infty})$ is well defined.\
\end{definition}

\begin{definition}\label{MixingAndRigidSequencesDefinition}
    Let $(x_n)_{n = 1}^\infty$ be a bounded sequence of complex numbers.
    \begin{enumerate}[(i)]
        \item If for any permissible triple $((x_n)_{n = 1}^\infty,$ $(y_n)_{n = 1}^\infty,(N_q)_{q = 1}^\infty)$ we have

        \begin{equation}\label{NearlyWeaklyMixingConditionEquation}
            \lim_{H\rightarrow\infty}\frac{1}{H}\sum_{h = 1}^H\lim_{q\rightarrow\infty}\left|\frac{1}{N_q}\sum_{n = 1}^{N_q}x_{n+h}\overline{y_n}\right| = 0,
        \end{equation}
        then $(x_n)_{n = 1}^\infty$ is a \textbf{nearly weakly mixing sequence}.

        \item If for any permissible triple $((x_n)_{n = 1}^\infty,(x_n)_{n = 1}^\infty,(N_q)_{q = 1}^\infty)$ there is a function $f \in L^1([0,1],m)$ for which

        \begin{equation}\label{SpectrallyLebesgueConditionEquation}
            \lim_{q\rightarrow\infty}\left|\frac{1}{N_q}\sum_{n = 1}^{N_q}x_{n+h}\overline{x_n}\right| = \int_0^1f(x)e^{2\pi ihx}dx,
        \end{equation}
        then $(x_n)_{n = 1}^\infty$ is a \textbf{spectrally Lebesgue sequence}.

        \item If for any permissible triple $((x_n)_{n = 1}^\infty,(x_n)_{n = 1}^\infty,(N_q)_{q = 1}^\infty)$ we have

        \begin{equation}
            \lim_{q\rightarrow\infty}\left|\frac{1}{N_q}\sum_{n = 1}^{N_q}x_{n+h}\overline{x_n}\right| = 0,
        \end{equation}
        then $(x_n)_{n = 1}^\infty$ is a \textbf{nearly orthogonal sequence}.
        \item If for any permissible triple of the form $((x_n)_{n = 1}^\infty,(x_n)_{n = 1}^\infty,(N_q)_{q = 1}^\infty)$ we have

        \begin{equation}
            \underset{m \in \mathbb{N}}{\sup }\ \underset{1 \le k \le K}{\text{min }}\ \lim_{q\rightarrow\infty}\frac{1}{N_q}\sum_{n = 1}^{N_q}|c_{n+m}-c_{n+k}|^2 < \epsilon.
        \end{equation}
        then $(x_n)_{n = 1}^\infty$ is a \textbf{compact sequence}\footnote{This definition was motivated by Definition 3.13 in \cite{TheErdosSumsetPaper}.}.
        \item If for any permissible triple of the form $((x_n)_{n = 1}^\infty,(x_n)_{n = 1}^\infty,(N_q)_{q = 1}^\infty)$ we have

        \begin{equation}
            \lim_{q\rightarrow\infty}\left|\frac{1}{N_q}\sum_{n = 1}^{N_q}x_{n+h}\overline{x_n}\right| = \int_0^1e^{2\pi ihx}d\mu(x),
        \end{equation}
        where $\mu\perp m$ is a positive measure, then $(x_n)_{n = 1}^\infty$ is a \textbf{spectrally singular sequence}.
    \end{enumerate}
\end{definition}

It is worth observing that if $(x_n)_{n = 1}^\infty$ satisfies any of items $(i)-(v)$ in Definition \ref{MixingAndRigidSequencesDefinition}, then $(\overline{x}_n)_{n = 1}^\infty$ satisfies the same list of items.

\begin{definition}\label{MixingAndRigidDistributionsDefinitions}
    Let $d \in \mathbb{N}, \mathcal{C}_d = \{f \in C([0,1]^d)\ |\ \int_{[0,1]^d}fdm^d = 0\}$, and $(x_n)_{n = 1}^\infty \subseteq [0,1]^d$.
    \begin{enumerate}[(i)]
        \item If for every $f \in \mathcal{C}_d$, $(f(x_n)_{n = 1}^\infty)$ is a nearly weakly mixing sequence, then $(x_n)_{n = 1}^\infty$ is a \textbf{wm-sequence}.
        \item If for every $f \in \mathcal{C}_d$, $(f(x_n)_{n = 1}^\infty)$ is a spectrally Lebesgue sequence, then $(x_n)_{n = 1}^\infty$ is a \textbf{sL-sequence}.
        \item If for every $f \in \mathcal{C}_d$, $(f(x_n)_{n = 1}^\infty)$ is a nearly orthogonal sequence, then $(x_n)_{n = 1}^\infty$ is a \textbf{o-sequence}.
        \item If for every $f \in C([0,1]^d)$, $(f(x_n)_{n = 1}^\infty)$ is a compact sequence, then $(x_n)_{n = 1}^\infty$ is a \textbf{c-sequence}.
        \item If for every $f \in C([0,1]^d)$, $(f(x_n)_{n = 1}^\infty)$ is a spectrally singular sequence, then $(x_n)_{n = 1}^\infty$ is a \textbf{ss-sequence}.
        
    \end{enumerate}
\end{definition}

\begin{definition}
\label{DefinitionOfNaturalDensity}
For a sequence of natural numbers $A = (n_k)_{k = 1}^{\infty}$ let

\begin{equation}
    \underline{d}(A) := \liminf_{N\rightarrow\infty}\frac{1}{N}\left|\{1 \le n \le N\ |\ n \in A\}\right| = \liminf_{N\rightarrow\infty}\frac{1}{N}\sum_{n = 1}^N\mathbbm{1}_A(n)\text{, and}
\end{equation}

\begin{equation}
    \overline{d}(A) := \limsup_{N\rightarrow\infty}\frac{1}{N}\left|\{1 \le n \le N\ |\ n \in A\}\right| = \limsup_{N\rightarrow\infty}\frac{1}{N}\sum_{n = 1}^N\mathbbm{1}_A(n).
\end{equation}
$\underline{d}(A)$ is the \textbf{natural lower density of A} and $\overline{d}(A)$ is the \textbf{natural upper density of A}. If $\overline{d}(A) = \underline{d}(A)$, then we let $d(A)$ denote the common value which is the \textbf{natural density of A}.
\end{definition}

\begin{definition} Let $A = (n_k)_{k = 1}^{\infty} \subseteq \mathbb{N}$ be a strictly increasing sequence satisfying $\underline{d}(A) > 0$. If $(\mathbbm{1}_A(n))_{n = 1}^\infty$ is a compact sequence as a sequence of bounded complex numbers, then $A$ is \textbf{compact} as a sequence of natural numbers. Similarly, if $(\mathbbm{1}_A(n))_{n = 1}^\infty$ is a spectrally singular as a sequence of bounded complex numbers, then $A$ is \textbf{spectrally singular} as a sequence of natural numbers.
\end{definition}

\subsubsection{A potpourri of lemmas}

Our first lemma follows from \cite[Theorem 2.1]{LebesgueSpectrumvanderCorput} and the preceding discussion.

\begin{lemma}\label{WhyWeNeedPermissibleTriples}
    Let $((x_n)_{n = 1}^\infty,$ $(y_n)_{n = 1}^\infty,(N_q)_{q = 1}^\infty)$ be a permissible triple. Then there exists a Hilbert space $\mathcal{H}$ and an injective map $\i$ from $\mathcal{H}$ to the space of sequence of complex numbers satisfying the following properties:
    \begin{enumerate}[(i)]
        \item For any $\xi_1,\xi_2 \in \mathcal{H}$ and $(z_{i,n})_{n = 1}^\infty := \i(\xi_i)$, we have

    \begin{equation}
        \langle \xi_1,\xi_2\rangle = \lim_{q\rightarrow\infty}\frac{1}{N_q}\sum_{n = 1}^{N_q}z_{1,n}\overline{z_{2,n}}.
    \end{equation}
    
        \item The map $\i$ is linear in the sense that for any $c \in \mathbb{C}$ and $(z_{3,n})_{n = 1}^\infty := \i(\xi_1+c\xi_2)$, we have

        \begin{equation}
            \lim_{q\rightarrow\infty}\frac{1}{N_q}\sum_{n = 1}^{N_q}|z_{3,n}-z_{1,n}-cz_{2,n}|^2 = 0.
        \end{equation}
        \item There exists a unitary operator $S:\mathcal{H}\rightarrow\mathcal{H}$ for which

    \begin{equation}
        \langle S^h\xi_1,\xi_2\rangle = \lim_{q\rightarrow\infty}\frac{1}{N_q}\sum_{n = 1}^{N_q}z_{1,n+h}\overline{z_{2,n}}.
    \end{equation}

        \item $\mathcal{H}$ is the smallest $S$-invariant Hilbert space containing the vectors $\vec{x}$ and $\vec{y}$ satisfying $\i(\vec{x}) = (x_n)_{n = 1}^\infty$ and $\i(\vec{y}) = (y_n)_{n = 1}^\infty$.
    \end{enumerate}
\end{lemma}

\begin{remark}
    Lemma \ref{WhyWeNeedPermissibleTriples} not only shows why we like to work with permissible triples, but it also justifies the terminology used in Definition \ref{MixingAndRigidSequencesDefinition}. To illustrate the latter claim with an example, let $(x_n)_{n = 1}^\infty$ be a spectrally singular sequence. We see that if $((x_n)_{n = 1}^\infty,(x_n)_{n = 1}^\infty,(N_q)_{q = 1}^\infty)$ is a permissible triple and $\vec{x}$ is as in Lemma \ref{WhyWeNeedPermissibleTriples}(iv), then

    \begin{equation}
        \langle S^h\vec{x},\vec{x}\rangle = \lim_{q\rightarrow\infty}\frac{1}{N_q}\sum_{n = 1}^{N_q}x_{n+h}\overline{x_n},
    \end{equation}
    hence the spectral measure of $\vec{x}$ with respect to $S$ is singular with respect to $m$.
\end{remark}

Our next lemma is a simple extension of \cite[Lemma 3.26]{TheErdosSumsetPaper}.

\begin{lemma}\label{RepresentationLemma}
Let $(c_n)_{n = 1}^{\infty} \subseteq \mathbb{C}$ be bounded. There exists a compact metric space $Y$, a continuous map $S:X\rightarrow X$, a continuous function $F:X\rightarrow\mathbb{C}$, and a point $x \in X$ with a dense orbit under $S$ such that $c_n = F(S^nx)$ for all $n \in \mathbb{N}$. Furthermore, if $(c_n)_{n = 1}^\infty$ is compact, then $(X,S)$ is such that for any weak$^*$ limit $\mu$ of $\{\frac{1}{N}\sum_{n = 1}^N\delta_{S^nx}\}_{N = 1}^\infty$, the measure preserving system $(X,S,\mu)$ has discrete spectrum.
\end{lemma}

\begin{proof}
    We may assume without loss of generality that $|c_n| \le 1$ for all $n \in \mathbb{N}$. Let $B \subseteq \mathbb{C}$ denote the closed unit ball and let $X' = B^\mathbb{N}$ endowed with the product topology. Since $B$ is a compact metric space we see that $X'$ is also a compact metric space. Let $S:X'\rightarrow X'$ denote the left shift. Let $F:X'\rightarrow\mathbb(C)$ denote the projection onto the first coordinate, which is seen to be continuous. Let $x = (c_n)_{n = 1}^\infty \in X'$, and let $X = c\ell(\{S^nx\}_{n = 1}^\infty)$. Since $X$ is a closed subset of $X'$ we see that $X$ is also a compact metric space, and it is clear that $x$ has a dense orbit in $X$ by construction. Since $F(S^nx) = c_n$, the first part of the lemma is proven. To see the latter claim, let $(N_q)_{q = 1}^\infty$ be such that

    \begin{equation}
        \lim_{q\rightarrow\infty}\frac{1}{N_q}\sum_{n = 1}^{N_q}\delta_{S^nx} = \mu,
    \end{equation}
    with convergence taking place in the weak$^*$ topology. We see that for any $h \in \mathbb{N}$, we have

    \begin{equation}
        \lim_{q\rightarrow\infty}\frac{1}{N_q}\sum_{n = 1}^{N_q}c_{n+h}\overline{c_n} = \int_XS^hF\cdot Fd\mu,
    \end{equation}
    so $((c_n)_{n = 1}^\infty,(c_n)_{n = 1}^\infty,(N_q)_{q = 1}^\infty)$ is a permissible triple. It follows that for any $\epsilon > 0$, there exists $K \in \mathbb{N}$ such that

    \begin{alignat*}{2}
        &\sup_{m \in \mathbb{N}}\min_{1 \le k \le K}||S^mF-S^kF||_2 = \sup_{m \in \mathbb{N}}\min_{1 \le k \le K}\int_X\left|S^mF-S^kF\right|^2d\mu\\
        =&\underset{m \in \mathbb{N}}{\sup }\ \underset{1 \le k \le K}{\text{min }}\ \lim_{q\rightarrow\infty}\frac{1}{N_q}\sum_{n = 1}^{N_q}|F(S^{n+m}x)-F(S^{n+k}x)|^2\\
        =&\underset{m \in \mathbb{N}}{\sup }\ \underset{1 \le k \le K}{\text{min }}\ \lim_{q\rightarrow\infty}\frac{1}{N_q}\sum_{n = 1}^{N_q}|c_{n+m}-c_{n+k}|^2 < \epsilon,
    \end{alignat*}
    so $F$ is a compact vector. Since the set of bounded compact vectors is a norm-closed $S$ invariant algebra, we see that $C(X)$ consists of compact vectors, hence all of $L^2(X,\mu)$ consists of compact vectors, which yields the desired result.
\end{proof}

\begin{lemma}\label{CreatingSpectrallyLebesgue/NWMSequences}
    Let $(x_n)_{n = 1}^\infty$ be a bounded sequence of complex numbers.
    \begin{enumerate}[(i)]
        \item If 

        \begin{equation}
            \sum_{h = 1}^\infty\limsup_{N\rightarrow\infty}\left|\frac{1}{N}\sum_{n = 1}^Nx_{n+h}\overline{x_n}\right|^2 < \infty,
        \end{equation}
        then $(x_n)_{n = 1}^\infty$ is spectrally Lebesgue.

        \item If 

        \begin{equation}
            \lim_{H\rightarrow\infty}\sum_{h = 1}^H\limsup_{N\rightarrow\infty}\left|\frac{1}{N}\sum_{n = 1}^Nx_{n+h}\overline{x_n}\right| = 0,
        \end{equation}
        then $(x_n)_{n = 1}^\infty$ is a nearly weakly mixing sequence.
    \end{enumerate}
\end{lemma}

\begin{proof}[Proof of (i)]
    Item (i) is a special case of \cite[Theorem 2.7]{LebesgueSpectrumvanderCorput} and item (ii) is a special case of \cite[Corollary 2.2.11]{SohailsPhDThesis}.
\end{proof}

\begin{lemma}
    Let $\mathcal{C}$ denote one of the following classes of bounded sequences of complex numbers: nearly weakly mixing, spectrally Lebesgue, compact, spectrally singular.
    \begin{enumerate}[(i)]
        \item If $(x_n)_{n = 1}^\infty,$ $(y_n)_{n = 1}^\infty \in \mathcal{C}$ and $c \in \mathbb{C}$, then for $z_n = x_n+cy_n$ we have that $(z_n)_{n = 1}^\infty \in \mathcal{C}$.

        \item If $\{(z_{n,k})_{n = 1}^\infty\}_{k = 1}^\infty \subseteq \mathcal{C}$ and there exists a bounded sequence $(z_n)_{n = 1}^\infty$ satisfying

        \begin{equation}
            \lim_{k\rightarrow\infty}\sup_{n \in \mathbb{N}}|z_{n,k}-z_n| = 0,
        \end{equation}
        then $(z_n)_{n = 1}^\infty \in \mathcal{C}$.
    \end{enumerate}
\end{lemma}

\begin{proof}
    We first prove (i). Let $((z_n)_{n = 1}^\infty,(z_n)_{n = 1}^\infty,(N_q')_{q = 1}^\infty)$ be a permissible triple, let $(N_q)_{q = 1}^\infty$ be a subsequence of $(N_q')_{q = 1}^\infty$ for which $((x_n)_{n = 1}^\infty,$ $(y_n)_{n = 1}^\infty,$ $(N_q)_{q = 1}^\infty)$ is a permissible triple, and let $\mathcal{H},\vec{x},\vec{y},$ and $S$ be as in Lemma \ref{WhyWeNeedPermissibleTriples}. We see that the spectral measures of $\vec{x}$ and $\vec{y}$ are continuous if $C$ denotes nearly weakly mixing sequences, absolutely continuous to Lebesgue if $\mathcal{C}$ denotes spectrally Lebesgue sequences, discrete if $\mathcal{C}$ denotes compact sequences, and mutually singular with Lebesgue if $\mathcal{C}$ denotes spectrally singular sequences. Consequently, the spectral measure of $\vec{x}+c\vec{y}$ is of the same class as that of $\vec{x}$ (and $\vec{y}$). 

    We see that (ii) is immediate if $\mathcal{C}$ denotes either the class of nearly weakly mixing sequences, or the class of compact sequences. Now we prove (ii) by way of contradiction when $\mathcal{C}$ denotes the class of spectrally Lebesgue sequences. Since $(z_n)_{n = 1}^\infty$ is not spectrally Lebesgue, let $((z_n)_{n = 1}^\infty,(z_n)_{n = 1}^\infty,(N_q')_{q = 1}^\infty)$ be a permissible triple for which

    \begin{equation}
        \gamma(h) := \lim_{q\rightarrow\infty}\frac{1}{N_q'}\sum_{n = 1}^{N_q'}z_{n+h}\overline{z_n}
    \end{equation}
    form the Fourier coefficients of some measure $\mu$ that is not absolutely continuous with respect to $m$. Let us apply Lemma \ref{WhyWeNeedPermissibleTriples} to the permissible triple $((z_n)_{n = 1}^\infty,$ $(z_n)_{n = 1}^\infty,$ $(N_q')_{q = 1}^\infty)$ to obtain $\mathcal{H}'$, $\vec{z}$, and $S'$. The spectral theorem lets us write $\vec{z} = \vec{z}_L+\vec{z}_s$ where $\nu_{\vec{z}_L} << m$ and $0 \neq \nu_{\vec{z}_s}\perp m$. Let $0 < \epsilon < ||\vec{z}_s||$ be arbitrary, let $k \in \mathbb{N}$ be such that $\sup_{n \in \mathbb{N}}|z_{n,k}-z_n| < \epsilon$, and let $(N_q)_{q = 1}^\infty$ be a subsequence of $(N_q')_{q = 1}^\infty$ for which $((z_n)_{n = 1}^\infty,(z_{n,k})_{n = 1}^\infty,(N_q)_{q = 1}^\infty)$ is a permissible triple. Let $\mathcal{H},\vec{z},\vec{z}_k,$ and $S$ be given by Lemma \ref{WhyWeNeedPermissibleTriples}. We see that $\mathcal{H}'$ is naturally identified with a subspace of $\mathcal{H}$ and $S'$ is the restriction of $S$ to $\mathcal{H}'$, which is justifies our abuse of notation for the term $\vec{z}$ in both situations. Since $\nu_{\vec{z}_k} << m$, we see that

    \begin{equation}
        \epsilon < ||\vec{z}_s|| \le ||\vec{z}-\vec{z}_k|| < \epsilon,
    \end{equation}
    which yields the desired contradiction. The proof in the case that $\mathcal{C}$ denotes the class of spectrally singular sequences is similar.
\end{proof}

\begin{corollary}\label{CheckingForwm/sL/c/ss}
    Let $d \in \mathbb{N}$ and $(x_n)_{n = 1}^\infty \subseteq [0,1]^d$.
    \begin{enumerate}[(i)]
        \item If for every $\vec{m} \in \mathbb{Z}^d\setminus\{0\}$ the sequence $(e(\vec{m}\cdot x_n))_{n = 1}^\infty$ is nearly weakly mixing, then $(x_n)_{n = 1}^\infty$ is a wm-sequence.
        \item If for every $\vec{m} \in \mathbb{Z}^d\setminus\{0\}$ the sequence $(e(\vec{m}\cdot x_n))_{n = 1}^\infty$ is spectrally Lebesgue, then $(x_n)_{n = 1}^\infty$ is a sL-sequence.
        \item If for every $\vec{m} \in \mathbb{Z}^d\setminus\{0\}$ the sequence $(e(\vec{m}\cdot x_n))_{n = 1}^\infty$ is compact, then $(x_n)_{n = 1}^\infty$ is a c-sequence.
        \item If for every $\vec{m} \in \mathbb{Z}^d\setminus\{0\}$ the sequence $(e(\vec{m}\cdot x_n))_{n = 1}^\infty$ is spectrally singular, then $(x_n)_{n = 1}^\infty$ is a ss-sequence.
    \end{enumerate}
\end{corollary}

\begin{lemma}\label{ErgodicDichotomyForSequences}
    Let $(x_n)_{n = 1}^\infty$ and $(y_n)_{n = 1}^\infty$ be bounded sequences of complex numbers.
    \begin{enumerate}[(i)]
        \item\label{LebesguePerpSingular} If $(x_n)_{n = 1}^\infty$ is spectrally Lebesgue and $(y_n)_{n = 1}^\infty$ is spectrally singular, then

        \begin{equation}
            \lim_{N\rightarrow\infty}\frac{1}{N}\sum_{n = 1}^Nx_ny_n = 0.
        \end{equation}
        \item\label{WeakMixingPerpCompact} If $(x_n)_{n = 1}^\infty$ is nearly weakly mixing and $(y_n)_{n = 1}^\infty$ is compact, then

        \begin{equation}
            \lim_{N\rightarrow\infty}\frac{1}{N}\sum_{n = 1}^Nx_ny_n = 0.
        \end{equation}
    \end{enumerate}
\end{lemma}

\begin{proof}
    Let $(N_q')_{q = 1}^\infty \subseteq \mathbb{N}$ be any sequence for which
    \begin{equation}
        I := \lim_{q\rightarrow\infty}\frac{1}{N_q'}\sum_{n = 1}^{N_q'}x_ny_n = 0
    \end{equation}
    exists. Let $(N_q)_{q = 1}^\infty$ be any subsequence of $(N_q')_{q = 1}^\infty$ for which $((x_n)_{n = 1}^\infty,$ $(y_n)_{n = 1}^\infty,$ $(N_q)_{q = 1}^\infty)$ is a permissible triple. Let $\mathcal{H},\vec{x},\vec{y}$, and $S$ be as in Lemma \ref{WhyWeNeedPermissibleTriples}. As mentioned in Section \ref{SpectralTheorySubsection}, we see that $\nu_{\vec{x}}\perp\nu_{\vec{y}}$, so $I = \langle \vec{x},\vec{y}\rangle = 0$.
\end{proof}

\begin{lemma}\label{MakingCompactSequences}
    Let $(c_n)_{n = 1}^\infty$ be a bounded and compact sequence of complex numbers.
    \begin{enumerate}[(i)]
        \item For all $\ell \in \mathbb{N}$, $(c_n^\ell)_{n = 1}^\infty$ is a compact sequence.

        \item If $(c_n)_{n = 1}^\infty$ takes nonnegative values, then $(\sqrt{c_n})_{n = 1}^\infty$ is a compact sequence.

        \item The sequence $(c_n')_{n = 1}^\infty$ given by $c_n'\overline{c_n} = |c_n|$ is a compact sequence.
    \end{enumerate}
\end{lemma}

\begin{proof}
    In all cases we may assume without loss of generality that $|c_n| \le 1$ for all $n$. To prove (i), it suffices to observe that $|a^\ell-b^\ell| = |a-b|\cdot|a^{\ell-1}+a^{\ell-2}b+\cdots+ab^{\ell-2}+b^{\ell-1}| \le \ell\max(|a|,|b|)^{\ell-1}|a-b|$, so 

    \begin{equation}
        \lim_{q\rightarrow\infty}\frac{1}{N_q}\sum_{n = 1}^{N_q}|c_{n+m}^\ell-c_{n+k}^\ell|^2 \le \ell^2\lim_{q\rightarrow\infty}\frac{1}{N_q}\sum_{n = 1}^{N_q}|c_{n+m}-c_{n+k}|^2.
    \end{equation}

    To prove (ii), let $\epsilon > 0$ be arbitrary, let $((c_n)_{n = 1}^\infty,(c_n)_{n = 1}^\infty,(N_q)_{q = 1}^\infty)$ be a permissible trip be a permissible triple, and let $K \in \mathbb{N}$ be such that

    \begin{equation}
        \sup_{m \in \mathbb{N}}\min_{1 \le k \le K}\lim_{q\rightarrow\infty}\frac{1}{N_q}\sum_{n = 1}^{N_q}|c_{n+m}-c_{n+k}|^2 < \epsilon^2.
    \end{equation}
    For $m \in \mathbb{N}$, let $k = k(m) \in [1,K]$ attain the minimum, and let $A = \{n \in \mathbb{N}\ |\ |\sqrt{c_{n+m}}+\sqrt{c_{n+k}}| \ge \sqrt{\epsilon}\}$. We see that 

    \begin{alignat*}{2}
        &\lim_{q\rightarrow\infty}\frac{1}{N_q}\sum_{n = 1}^{N_q}|\sqrt{c_{n+m}}-\sqrt{c_{n+k}}|^2\\
        \le &\lim_{q\rightarrow\infty}\frac{1}{N_q}\left(\sum_{n \in [1,N_q]\cap A}|\sqrt{c_{n+m}}-\sqrt{c_{n+k}}|^2+\sum_{n \in [1,N_q]\cap A^c}|\sqrt{c_{n+m}}-\sqrt{c_{n+k}}|^2\right)\\
        \le &\lim_{q\rightarrow\infty}\frac{1}{N_q}\left(\sum_{n \in [1,N_q]\cap A}\left|\frac{c_{n+m}-c_{n+k}}{\sqrt{c_{n+m}}+\sqrt{c_{n+k}}}\right|^2+\sum_{n \in [1,N_q]\cap A^c}|\sqrt{c_{n+m}}+\sqrt{c_{n+k}}|^2\right)\\
        \le &\frac{1}{\epsilon}\lim_{q\rightarrow\infty}\frac{1}{N_q}\sum_{n = 1}^{N_q}|c_{n+m}-c_{n+k}|^2+\epsilon \le 2\epsilon.
    \end{alignat*}

    To prove (iii), we recall that $(\overline{c_n})_{n = 1}^\infty$ is a compact sequence, so $(|c_n|^2 = c_n\overline{c_n})_{n = 1}^\infty$ is a compact sequence by \cite[Lemma 2.3.8(ii)]{SohailsPhDThesis}, so the desrired result now follows from Lemma \ref{MakingCompactSequences}(ii).
\end{proof}
\section{Main Results}\label{MainResultsSection}

\begin{theorem}\label{PropertiesOfsLSequences}
    Let $d,d_1 \in \mathbb{N}$ and let $(x_n)_{n = 1}^\infty \subseteq [0,1]^d$. 
    \begin{enumerate}[(i)]
        \item If $(x_n)_{n = 1}^\infty$ is a sL-sequence (wm-sequence) and $(n_k)_{k = 1}^\infty \subseteq \mathbb{N}$ is a spectrally singular (compact) sequence, then $(x_{n_k})_{k = 1}^\infty$ is uniformly distributed.
        
        \item If for any compact sequence $(n_k)_{k = 1}^\infty \subseteq \mathbb{N}$, the sequence $(x_{n_k})_{k = 1}^\infty$ is uniformly distributed, then $(x_n)_{n = 1}^\infty$ is a wm-sequence.

        \item If $(x_n)_{n = 1}^\infty$ is a sL-sequence (wm-sequence) and $(y_n)_{n = 1}^\infty \subseteq [0,1]^{d_1}$ is a uniformly distributed ss-sequence (c-sequence), then $(x_n,y_n)_{n = 1}^\infty$ is uniformly distributed.

        \item If $(x_n+y_n)_{n = 1}^\infty$ is uniformly distributed for all c-sequences $(y_n)_{n = 1}^\infty \subseteq [0,1]^d$, then $(x_n)_{n = 1}^\infty$ is a wm-sequence.
    \end{enumerate}
\end{theorem}

\begin{proof}[Proof of (i)]
Let $\vec{m} \in \mathbb{Z}^d\setminus\{\vec{0}\}$ be arbitrary, and note that $(e^{2\pi i\vec{m}\cdot x_n})_{n = 1}^{\infty}$ is a spectrally Lebesgue (nearly weakly mixing) sequence . Letting $B = (n_k)_{k = 1}^{\infty}$, we see that $(\mathbbm{1}_B(n))_{n = 1}^{\infty}$ is a spectrally singular sequence (compact sequence), so by Lemma \ref{ErgodicDichotomyForSequences} we see that

\begin{alignat*}{2}
    0 & = \lim_{N\rightarrow\infty}\frac{1}{N}\left|\sum_{n = 1}^Ne(\vec{m}\cdot x_n)\mathbbm{1}_B(n)\right| = \lim_{N\rightarrow\infty}\frac{1}{N}\left|\sum_{n_k \in [1,N]}e(\vec{m}\cdot x_{n_k})\right|\\
    &\ge \lim_{K\rightarrow\infty}\frac{\underline{d}(B)}{K}\left|\sum_{k = 1}^Ke(\vec{m}\cdot x_{n_k})\right|.
\end{alignat*}
\end{proof}

\begin{proof}[Proof of (ii)]
Let $\vec{m} = (m_1,\cdots,m_d) \in \mathbb{Z}^d\setminus\{0\}$ be arbitrary, let $((e(\vec{m}\cdot x_n))_{n = 1}^\infty,(e(\vec{m}\cdot x_n))_{n = 1}^\infty,(N_q)_{q = 1}^\infty)$ be a permissible triple, and let $\mathcal{H}$, $S$, and $\i$ be as in Lemma \ref{WhyWeNeedPermissibleTriples}. Let $\xi := \i^{-1}((e(\vec{m}\cdot x_n))_{n = 1}^\infty)$, and let $\xi = \xi_w+\xi_c$ where $\xi_c$ has a pre-compact orbit under $S$ in the norm topology, and $\xi_w$ is weakly mixing with respect to $S$. Letting $(w_n)_{n = 1}^\infty = \i(\xi_w)$ and $(C_n)_{n = 1}^\infty = \i(\xi_c)$, we can assume without loss of generality that $e(\vec{m}\cdot x_n) = C_n+d_n$. Using \cite[Lemma 2.2.2]{SohailsPhDThesis} and \cite[Theorem 2.25]{ProjectionLemma}, we see that $\sup_n|C_n| \le \sup_n|e(\vec{m}\cdot x_n)| = 1$. Let $(c_n)_{n = 1}^\infty \subseteq \mathbb{S}^1$ be such that $c_n\overline{C_n} = |C_n|$. We see that $(c_n)_{n = 1}^\infty$ is a compact sequence by Lemma \ref{MakingCompactSequences}(iii). Let $x,X,S,$ and $F$ be as in Lemma \ref{RepresentationLemma} and satisfy $F(S^nx) = c_n$. Now let $(N_q)_{q = 1}^\infty$ be such that

\begin{equation}
    \lim_{q\rightarrow\infty}\frac{1}{N_q}\sum_{n = 1}^{N_q}\delta_{S^nx} = \mu,
\end{equation}
and recall that $(X,S,\mu)$ has discrete spectrum by Lemma \ref{RepresentationLemma}. Now let $\epsilon > 0$ be arbitrary, and let $g := \sum_{i = 1}^Md_i\mathbbm{1}_{A_i}$ be such that $||F-g||_\infty < \epsilon$ and each $A_i$ is open. For $1 \le i \le M$, let $(c_{i,n})_{n = 1}^\infty$ be given by $c_{i,n} = \mathbbm{1}_{A_i}(S^nx)$. To see that $(c_{i,n})_{n = 1}^\infty$ is compact, let $((c_{i,n})_{n = 1}^\infty,(c_{i,n})_{n = 1}^\infty,(N_q')_{q = 1}^\infty)$ be a permissible triple, and let $(N_q'')_{q = 1}^\infty$ be any subsequence of $(N_q')_{q = 1}^\infty$ for which

\begin{equation}
    \lim_{q\rightarrow\infty}\frac{1}{N_q''}\sum_{n = 1}^{N_q''}\delta_{S^nx} = \mu',
\end{equation}
with convergence taking place in the weak topology. Lemma \ref{RepresentationLemma} tells us that $(X,S,\mu')$ has discrete spectrum, so $\mathbbm{1}_{A_i}$ has precompact orbit under $S$, which yields the desired result. It follows that if we let $(n_{i,k})_{k = 1}^\infty$ be an increasing enumeration of those $(c_{i,k})_{k = 1}^\infty$ that are $1$, then each $(n_{i,k})_{k = 1}^\infty$ would be a compact sequence of natural numbers if they have positive lower density, which is not necessarily the case. To overcome this difficulty, we let $(n_{i,k,1})_{k = 1}^\infty$ be the increasing enumeration of $(n_{i,k})_{k = 1}^\infty\cup2\mathbb{N}$ and $(n_{i,k,2})_{k = 1}^\infty$ be the increasing enumeration of $(n_{i,k})_{k = 1}^\infty\cup(2\mathbb{N}+1)$, so that $(n_{i,k,j})_{k = 1}^\infty$ is a compact sequence of natural numbers for $j = 1,2$. Since $(x_{n_k})_{k = 1}^\infty$ is uniformly distributed for any compact sequence $(n_k)_{k = 1}^\infty$, we see that for $A_i = (n_{i,k})_{k = 1}^\infty$ and $A_{i,j} := (n_{i,k,j})_{k = 1}^\infty$ we have

\begin{alignat*}{2}
    0&=\lim_{N\rightarrow\infty}\frac{1}{N}\sum_{n = 1}^Ne(\vec{m}\cdot x_n)\mathbbm{1}_{A_{i,j}}(n) = \lim_{N\rightarrow\infty}\frac{1}{N}\sum_{n = 1}^Ne(\vec{m}\cdot x_n)\text{, hence}\\
    0&=\lim_{N\rightarrow\infty}\frac{1}{N}\sum_{n = 1}^Ne(\vec{m}\cdot x_n)\left(\mathbbm{1}_{A_{i,1}}(n)+\mathbbm{1}_{A_{i,2}}(n)-1\right)\\
    &=\lim_{N\rightarrow\infty}\frac{1}{N}\sum_{n = 1}^Ne(\vec{m}\cdot x_n)\mathbbm{1}_{A_{i}}(n)\text{, hence}\\
    0&= \lim_{N\rightarrow\infty}\frac{1}{N}\sum_{n = 1}^Ne(\vec{m}\cdot x_n)\sum_{i = 1}^md_i\mathbbm{1}_{A_i}(n).
\end{alignat*}
After recalling that $\sum_{i = 1}^md_i\mathbbm{1}_{A_i}(n) = g(S^nx)$, $|g(S^nx)-F(S^nx)| < \epsilon$, and that $\epsilon > 0$ was arbitrary, we see that

\begin{alignat*}{2}
    0 &= \lim_{N\rightarrow\infty}\frac{1}{N}\sum_{n = 1}^Ne(\vec{m}\cdot x_n)c_n = \lim_{N\rightarrow\infty}\frac{1}{N}\sum_{n = 1}^N(C_nc_n+d_nc_n)\\
    &=\lim_{N\rightarrow\infty}\frac{1}{N}\sum_{n = 1}^N|C_n| \ge \limsup_{N\rightarrow\infty}\frac{1}{N}\sum_{n = 1}^N|C_n|^2,
\end{alignat*}
so we may assume without loss of generality that $C_n = 0$ for all $n$, so $(e(\vec{m}\cdot x_n))_{n = 1}^\infty$ is a nearly weakly mixing sequence.
\end{proof}

\begin{proof}[Proof of (iii)]
Let $(\vec{m},\vec{m}_1) \in \mathbb{Z}^{d+d_1}\setminus\{\vec{0}\}$ be arbitrary. Note that $(e(\vec{m}\cdot x_n))_{n = 1}^\infty$ is a spectrally Lebesgue (nearly weakly mixing) sequence when $\vec{m} \neq \vec{0}$, while $(e(\vec{m}_1\cdot y_n))_{n = 1}^\infty$ is a spectrally singular (compact) sequence. Since $(y_n)_{n = 1}^\infty$ is uniformly distributed, it is clear that if $\vec{m} = \vec{0}$, then $\vec{m}_1 \neq \vec{0}$ and

\begin{equation}
    \lim_{N\rightarrow\infty}\frac{1}{N}\sum_{n = 1}^Ne((\vec{m}+\vec{m_1})\cdot(x_n,y_n)) = \lim_{N\rightarrow\infty}\frac{1}{N}\sum_{n = 1}^Ne(\vec{m_1}\cdot y_n) = 0.
\end{equation}
If $\vec{m} \neq \vec{0}$, then we can use Lemma \ref{ErgodicDichotomyForSequences} to see that

\begin{alignat*}{2}
    \lim_{N\rightarrow\infty}\frac{1}{N}\sum_{n = 1}^Ne((\vec{m}+\vec{m_1})\cdot(x_n,y_n)) = \lim_{N\rightarrow\infty}\frac{1}{N}\sum_{n = 1}^Ne(\vec{m}\cdot x_n)e(\vec{m}_1\cdot y_n) = 0.
\end{alignat*}
\end{proof}

\begin{proof}[Proof of (iv)]
Firstly, we see that $(0)_{n = 1}^\infty$ is a c-sequence, so $(x_n)_{n = 1}^\infty$ is uniformly distributed. Since we have already proven (ii), it suffices to show that $(x_{n_k})_{k = 1}^\infty$ is uniformly distributed for an arbitrary compact sequence $A = (n_k)_{k = 1}^\infty \subseteq \mathbb{N}$. For $t \in [0,1]^d$, let $(y_n(t))_{n = 1}^\infty$ be given by $y_n(t) = t\mathbbm{1}_A(t)$.  Since $(y_n(t))_{n = 1}^\infty$ is always a c-sequence, we see that $(x_n+y_n(t))_{n = 1}^\infty$ is uniformly distributed for any $t \in (0,1)$. It follows that for all $t \in [0,1]^d$ we have

\begin{alignat*}{2}
    &0 = \lim_{N\rightarrow\infty}\frac{1}{N}\sum_{n = 1}^Nf(x_n+y_n(t)) = \lim_{N\rightarrow\infty}\frac{1}{N}\sum_{n = 1}^Nf(x_n+t\mathbbm{1}_A(n))\text{, hence}\\
    &0 = \lim_{N\rightarrow\infty}\frac{1}{N}\sum_{n = 1}^Nf((x_n+t)-f(x_n))\mathbbm{1}_A(n).
\end{alignat*}
Consequently, if $\mu$ is any weak$^*$ limit point of $\{\frac{1}{N}\sum_{n = 1}^N\delta_{x_n}\mathbbm{1}_A(n)\}_{N = 1}^\infty$, then $\mu$ is translation invariant, hence $\mu$ is a constant multiple of the Lebesgue measure, which shows that $(x_{n_k})_{k = 1}^\infty$ is uniformly distributed.
\end{proof}

\begin{remark}
    It is natural to ask if items (ii) and (iv) of Theorem \ref{PropertiesOfsLSequences} have analogues for sL-sequences. While we do not have counterexamples showing that such analogues do not exist, we can explain why our current methods of proof do not extend to this situation. The proof of Theorem \ref{PropertiesOfsLSequences}(ii) relies on the fact that if $(X,\mathscr{B},\mu,T)$ is a measure preserving system, then there is a $\sigma$-algebra $\mathcal{K} \subseteq \mathscr{B}$ with respect to which all eigenfunctions\footnote{We recall that the span closure of eigenfunctions is precisely those functions whose orbit under $T$ is pre-compact in the norm topology of $L^2(X,\mu)$, which is why compact sequences can be associated to eigenfunctions.} of $T$ are measurable. However, there need not exists $\sigma$-algebra $\mathcal{S} \subseteq \mathscr{B}$ with respect to which all $f \in L^2(X,\mu)$ that have singular spectrum are measurable. Similarly, the proof of Theorem \ref{PropertiesOfsLSequences}(iv) (intuitively) uses the fact that if $f$ is in the span-closure of eigenfunctions of $T$, then so it $f^k$ for any $k \in \mathbb{N}$, but the analogous statement for functions of singular spectrum is not true. 
\end{remark}

Our next two results can be seen as generalizations of Theorem \ref{UniformDistributionvanderCorput'sDifferenceTheorem}.

\begin{theorem}\label{GeneralizationOfClassicalvdC}
Let $d\in \mathbb{N}$. If $(x_n)_{n = 1}^\infty \subseteq [0,1]^d$ is such that 

\begin{equation}
    \sum_{h = 1}^\infty\overline{D}((x_{n+h}-x_n)_{n = 1}^\infty)^2 < \infty,
\end{equation}
then $(x_n)_{n = 1}^\infty$ is a sL-sequence.
\end{theorem}

\begin{proof}
    Let $\vec{m} = (m_1,\cdots,m_d) \in \mathbb{Z}^d$ be arbitrary and let $M = \max_{1 \le i \le d}(|m_i|)$. Using Theorem \ref{SpecialCaseOfKHInequality}, we see that for $h \in \mathbb{N}$ we have

\begin{alignat*}{2}
    &\left|\limsup_{N\rightarrow\infty}\frac{1}{N}\sum_{n = 1}^Ne(\vec{m}\cdot(x_{n+h}-x_n))\right| \le (4\pi M)^d\overline{D}\left((x_{n+h}-x_n)_{n = 1}^\infty\right)\text{, hence}\\
    &\sum_{h = 1}^\infty\left|\limsup_{N\rightarrow\infty}\frac{1}{N}\sum_{n = 1}^Ne(\vec{m}\cdot(x_{n+h}-x_n))\right|^2\\
    \le &(4\pi M)^{2d}\sum_{h = 1}^\infty\overline{D}\left((x_{n+h}-x_n)_{n = 1}^\infty\right)^2 < \infty,
\end{alignat*}
so Lemma \ref{CreatingSpectrallyLebesgue/NWMSequences}(i) tells us that $(e(\vec{m}\cdot x_n))_{n = 1}^\infty$ is a spectrally Lebesgue sequence. Since $\vec{m}$ was arbitrary, the desired result now follows from Corollary \ref{CheckingForwm/sL/c/ss}(ii).
\end{proof}

\begin{theorem}\label{WeakMixingUDvdC}
    Let $d \in \mathbb{N}$. If $(x_n)_{n = 1}^\infty \subseteq [0,1]^d$ is such that 

\begin{equation}
    \lim_{H\rightarrow\infty}\frac{1}{H}\sum_{h = 1}^H\overline{D}((x_{n+h}-x_n)_{n = 1}^\infty) = 0,
\end{equation}
then $(x_n)_{n = 1}^\infty$ is a wm-sequence.
\end{theorem}

\begin{proof}
    Let $\vec{m} = (m_1,\cdots,m_d) \in \mathbb{Z}^d$ be arbitrary and let $M = \max_{1 \le i \le d}(|m_i|)$. Using Theorem \ref{SpecialCaseOfKHInequality}, we see that for $h \in \mathbb{N}$ we have

\begin{alignat*}{2}
    &\left|\limsup_{N\rightarrow\infty}\frac{1}{N}\sum_{n = 1}^Ne(\vec{m}\cdot(x_{n+h}-x_n))\right| \le (4\pi M)^d\overline{D}\left((x_{n+h}-x_n)_{n = 1}^\infty\right)\text{, hence}\\
    &\lim_{H\rightarrow\infty}\frac{1}{H}\sum_{h = 1}^H\left|\limsup_{N\rightarrow\infty}\frac{1}{N}\sum_{n = 1}^Ne(\vec{m}\cdot(x_{n+h}-x_n))\right|\\
    \le &(4\pi M)^d\lim_{H\rightarrow\infty}\frac{1}{H}\sum_{h = 1}^H\overline{D}\left((x_{n+h}-x_n)_{n = 1}^\infty\right) = 0,
\end{alignat*}
so Lemma \ref{CreatingSpectrallyLebesgue/NWMSequences}(ii) tells us that $(e(\vec{m}\cdot x_n))_{n = 1}^\infty$ is a spectrally Lebesgue sequence. Since $\vec{m}$ was arbitrary, the desired result now follows from Corollary \ref{CheckingForwm/sL/c/ss}(i).
\end{proof}

\begin{remark}\label{RemarkJustifyingPermissiblePairs}
    In light of Lemma \ref{ExampleNecessitatingPermissibleTriples}, we see that Theorem \ref{GeneralizationOfClassicalvdC} does not characterize sL-sequences, and that Theorem \ref{WeakMixingUDvdC} does not characterize wm-sequences. A characterization of wm-sequences using permissible pairs has already been obtained in \cite{SohailsPhDThesis}, which we state below as Theorem \ref{4CharacterizationsOfwmSequences}. In light of Lemma \ref{ExampleNecessitatingPermissibleTriples}, we see, intuitively speaking, that the reason only reason Theorem \ref{WeakMixingUDvdC} is not a characterization of wm-sequences is that the $\limsup$ that is in the definition of $\overline{D}((x_{n+h}-x_n)_{n = 1}^{\infty})$ is attained by a sequence $(N_q(h))_{q = 1}^\infty$ that depends on $h$. This also indicates that the utility of permissible pairs is that they do not allow the sequence $(N_q)_{q = 1}^\infty$ to change with $h$.

    The situation for sL-sequences is similar, in that any potential characterization would have to use permissible pairs rather than $\overline{D}$. However, we like to have sufficient conditions using $\overline{D}$ rather than permissible pairs for the sake of application. Let us now give an example to show why we do not pursue a more general sufficient condition for sL-sequences than that present in Theorem \ref{GeneralizationOfClassicalvdC}. Recall that there exist positive measure $\nu\perp m$ on $[0,1]$ whose Fourier coefficients tend to 0. Consequently, we may pick a convex decreasing sequence of real numbers $(c_n)_{n = 1}^\infty$ satisfying $c_n \ge |\hat{\nu}(n)|$, and it is a classical result \cite[Theorem 4.1]{KatznelsonHarmonicAnalysis} that there exists a positive measure $\mu << m$ for which $\hat{\mu}(n) = c_n$. In particular, it is possible to decrease the Fourier coefficients of a measure $\mu << m$ and as a result obtain the Fourier coefficients of a measure $\nu\perp m$. However, this phenomenon does not happen if we start with a measure $\mu$ for which $\sum_{n = 1}^\infty|\hat{\mu}(n)|^2 < \infty.$
\end{remark}

\begin{theorem}[{\cite[Theorem 2.4.17]{SohailsPhDThesis}}]
\label{4CharacterizationsOfwmSequences}
For $(x_n)_{n = 1}^{\infty} \subseteq [0,1]^{d_1}$ the following are equivalent:

\begin{enumerate}[(i)]
    \item $(x_n)_{n = 1}^{\infty}$ is a wm-sequence.
    
    \item For all uniformly distributed $(y_n)_{n = 1}^{\infty} \subseteq [0,1]^{d_2}$ and $(N_q)_{q = 1}^{\infty} \subseteq \mathbb{N}$ for which \\ $(\{(x_n,y_{n+h})_{n = 1}^{\infty}\}_{h = 1}^{\infty},(N_q)_{q = 1}^{\infty})$ is a permissible pair, we have

\begin{equation}
    \lim_{H\rightarrow\infty}\frac{1}{H}\sum_{h = 1}^HD((x_n,y_{n+h})_{n = 1}^{\infty},(N_q)_{q = 1}^{\infty}) = 0.
    \label{EquationForWeakMixingDiscrepancy}
\end{equation}

\item For all $(N_q)_{q = 1}^{\infty} \subseteq \mathbb{N}$ for which $(\{(x_n,x_{n+h})_{n = 1}^{\infty}\}_{h = 1}^{\infty},(N_q)_{q = 1}^{\infty})$ is a permissible pair, we have

\begin{equation}
    \lim_{H\rightarrow\infty}\frac{1}{H}\sum_{h = 1}^HD((x_n,x_{n+h})_{n = 1}^{\infty},(N_q)_{q = 1}^{\infty}) = 0.
\end{equation}

\item For all $(N_q)_{q = 1}^{\infty} \subseteq \mathbb{N}$ for which $(\{(x_{n+h}-x_n)_{n = 1}^{\infty}\}_{h = 1}^{\infty},(N_q)_{q = 1}^{\infty}))$ is a permissible pair, we have

\begin{equation}
\label{WeakMixingDiscrepancy}
    \lim_{H\rightarrow\infty}\frac{1}{H}\sum_{h = 1}^HD((x_{n+h}-x_n)_{n = 1}^{\infty},(N_q)_{q = 1}^{\infty}) = 0.
\end{equation}
\end{enumerate}
\end{theorem}

In light of Remark \ref{RemarkJustifyingPermissiblePairs}, it is surprising that our next theorem characterizes o-sequences without the use of permissible pairs.

\begin{theorem}\label{CharacterizationOfOsequences} $(x_n)_{n = 1}^{\infty} \subseteq [0,1]^d$ is an $o$-sequence if and only if for each $h \in \mathbb{N}$,\\ $(x_n, x_{n+h})_{n = 1}^{\infty} \subseteq [0,1]^{2d}$ is uniformly distributed.
\end{theorem}

\begin{proof} For the first direction, let us assume that $(x_n,x_{n+h})_{n = 1}^{\infty}$ is uniformly distributed in $[0,1]^{2d}$ for all $h \in \mathbb{N}$. We see that for all $\vec{k}_1,\vec{k}_2 \in \mathbb{Z}^d$ that are not both $(0,0,\cdots,0)$ and any $h \in \mathbb{N}$ we have

\begin{equation}
    \lim_{N\rightarrow\infty}\frac{1}{N}\sum_{n = 1}^Ne(\vec{k_1}\cdot x_n+\vec{k}_2\cdot x_{n+h}) = 0.
\end{equation}
Now let $f \in C([0,1]^d)$ satisfy $\int_{[0,1]^d}fdm^d = 0$, let $\epsilon \in (0,||f||_\infty)$ be arbitrary, and let $K$ be such that

\begin{equation}
    ||f(x)-\sum_{k \in [-K,K]^d}c_{\vec{k}}e(\vec{k}\cdot x)||_{\infty} < \epsilon,
\end{equation}
for some $(c_{\vec{k}})_{\vec{k} \in [-K,K]^d}$. Since we may assume without loss of generality that $c_{(0,0,\cdots,0)} = 0$, we see that for all $h \in \mathbb{N}$ we have

\begin{alignat*}{2}
    & \lim_{N\rightarrow\infty}\left|\frac{1}{N}\sum_{n = 1}^Nf(x_{n+h})\overline{f(x_n)}\right|\le 3\epsilon||f||_{\infty}+\\
     & \qquad\lim_{N\rightarrow\infty}\left|\frac{1}{N}\sum_{n = 1}^N\left(\sum_{k \in [-K,K]^d}c_ke\left(\vec{k}\cdot x_{n+h}\right)\right)\left(\sum_{k \in [-K,K]^d}c_ke\left(-\vec{k}\cdot x_n\right)\right)\right|\\
    &= 3\epsilon||f||_{\infty}+\sum_{k_1,k_2 \in [-K,K]^d}\lim_{N\rightarrow\infty}\left|\frac{1}{N}\sum_{n = 1}^Nc_{k_1}\overline{c_{k_2}}e\left((\vec{k}_1,-\vec{k}_2)\cdot(x_{n+h},x_n)\right)\right|\\
    &= 3\epsilon||f||_{\infty}.
\end{alignat*}
Since $\epsilon > 0$ was arbitrary, we are done with the first direciton. For the reverse direction, let us assume that $(x_n)_{n = 1}^{\infty}$ is an o-sequence. We will first show that $(x_{n+h}-x_n)_{n = 1}^{\infty}$ is uniformly distributed for all $h \in \mathbb{N}$. To this end, let $\vec{k} \in \mathbb{Z}^d\setminus\{(0,0\cdots,0)\}$ and $h \in \mathbb{N}$ both be arbitrary and note that $(e(\vec{k}\cdot x_n))_{n = 1}^{\infty}$ is a nearly orthogonal sequence. Let $(N_q)_{q = 1}^{\infty}$ be any sequence for which

\begin{equation}
	\lim_{q\rightarrow\infty}\frac{1}{N_q}\sum_{n = 1}^{N_q}e\left(\vec{k}\cdot(x_{n+h}-x_n)\right)
\end{equation}
exists. By passing to a subsequence of $(N_q)_{q = 1}^{\infty}$ if necessary, we may assume without loss of generality that $\left(\left(e\left(\vec{k}\cdot x_n\right)\right)_{n = 1}^{\infty},\left(e\left(\vec{k}\cdot x_n\right)\right)_{n = 1}^{\infty},\left(N_q\right)_{q = 1}^{\infty}\right)$ is a permissible triple. Since $(e^{2\pi i\langle k, x_n\rangle})_{n = 1}^{\infty}$ is a nearly orthogonal sequence it follows that

\begin{equation}
	\lim_{q\rightarrow\infty}\frac{1}{N_q}\sum_{n = 1}^{N_q}e\left(\vec{k}\cdot(x_{n+h}-x_n)\right) = 0,
\end{equation}
from which it follows that $(x_{n+h}-x_n)_{n = 1}^{\infty}$ is indeed uniformly distributed for all $h \in \mathbb{N}$. Now let $h \in \mathbb{N}$ be arbitrary, let $\vec{k}_1,\vec{k}_2 \in \mathbb{Z}^d$ be such that $\vec{k}_1$ and $\vec{k}_2$ are not both $(0,0,\cdots,0)$ and let $(N_q)_{q = 1}^{\infty} \subseteq \mathbb{N}$ be such that

\begin{equation}\label{SomeEquation1}
	\lim_{q\rightarrow\infty}\frac{1}{N_q}\sum_{n = 1}^{N_q}e\left((\vec{k}_1,\vec{k}_2)\cdot(x_n,x_{n+h})\right)
\end{equation}
exists. If $\vec{k}_1$ or $\vec{k}_2$ is $(0,0,\cdots,0)$, then the limit in Equation \eqref{SomeEquation1} is $0$ since the o-sequence $(x_n)_{n = 1}^{\infty}$ is uniformly distributed, so let us assume that neither of $\vec{k}_1$ and $\vec{k}_2$ are $(0,0,\cdots,0)$. Note that for all $c \in \mathbb{C}$ we have that $(e(\vec{k}_1\cdot x_n)+ce(\vec{k}_2\cdot x_n))_{n = 1}^{\infty}$ is a nearly orthogonal sequence since $(x_n)_{n = 1}^\infty$ is an o-sequence, so we once again see that

\begin{alignat*}{2}
	0 & = \lim_{q\rightarrow\infty}\frac{1}{N_q}\sum_{q = 1}^{N_q}(e(\vec{k}_1\cdot x_{n+h})+ce(\vec{k}_2\cdot x_{n+h})(e(-\vec{k}_1\cdot x_n)+\overline{c}e(-\vec{k}_2\cdot x_n)\\
	& = \lim_{q\rightarrow\infty}\frac{1}{N_q}\sum_{q = 1}^{N_q}e(\vec{k}_1\cdot(x_{n+h}-x_n))+|c|^2\lim_{q\rightarrow\infty}\frac{1}{N_q}\sum_{q = 1}^{N_q}e(\vec{k}_2\cdot(x_{n+h}-x_n)) \\ &\textcolor{white}{=} +c\lim_{q\rightarrow\infty}\frac{1}{N_q}\sum_{q = 1}^{N_q}e(\vec{k}_2\cdot x_{n+h}-\vec{k}_1\cdot x_n)+\overline{c}\lim_{q\rightarrow\infty}\frac{1}{N_q}\sum_{q = 1}^{N_q}e(\vec{k}_1\cdot x_{n+h}-\vec{k}_2\cdot x_n)\\
	&= c\lim_{q\rightarrow\infty}\frac{1}{N_q}\sum_{q = 1}^{N_q}e(\vec{k}_2\cdot x_{n+h}-\vec{k}_1\cdot x_n)+\overline{c}\lim_{q\rightarrow\infty}\frac{1}{N_q}\sum_{q = 1}^{N_q}e(\vec{k}_1\cdot x_{n+h}-\vec{k}_2\cdot x_n).
\end{alignat*}

Letting $A(c)$ represent the final quantity in the previous calculation, we observe that

\begin{alignat*}{2}
	0 & = A(1)-iA(i)\\
	& = \lim_{q\rightarrow\infty}\frac{1}{N_q}\sum_{q = 1}^{N_q}e(\vec{k}_2\cdot x_{n+h}-\vec{k}_1\cdot x_n)+\lim_{q\rightarrow\infty}\frac{1}{N_q}\sum_{q = 1}^{N_q}e(\vec{k}_1\cdot x_{n+h}-\vec{k}_2\cdot x_n) \\ & \textcolor{white}{=} -i\Bigg(i\lim_{q\rightarrow\infty}\frac{1}{N_q}\sum_{q = 1}^{N_q}e(\vec{k}_2\cdot x_{n+h}-\vec{k}_1\cdot x_n)\\
 &\qquad\qquad\qquad-i\lim_{q\rightarrow\infty}\frac{1}{N_q}\sum_{q = 1}^{N_q}e(\vec{k}_1\cdot x_{n+h}-\vec{k}_2\cdot x_n)\Bigg)\\
	& = 2\lim_{q\rightarrow\infty}\frac{1}{N_q}\sum_{q = 1}^{N_q}e(\vec{k}_2\cdot x_{n+h}-\vec{k}_1\cdot x_n).
\end{alignat*}
Since $(k_1,k_2)\mapsto(-k_1,k_2)$ is a bijection from $\left(\mathbb{Z}^{d}\setminus\{(0,\cdots,0)\}\right)^2$ to itself, we see that $(x_n,x_{n+h})_{n = 1}^{\infty}$ is uniformly distributed for all $h \in \mathbb{N}$.
\end{proof}

It is worth noting that if $p(x) = a_nx^n+\cdots+a_1x+a_0 \in \mathbb{R}[x]$ is such that $n \ge 2$ and $a_j$ is irrational for some $j \ge 2$, then $(p(n))_{n = 1}^\infty$ is an o-sequence. In \cite[Theorem 2.4.30]{SohailsPhDThesis} is it shown that $(\lfloor g(n)\rfloor\alpha)_{n = 1}^\infty$ is an o-sequence for a large class of tempered functions $g$, such that $g(n) = n^t$ with $t \in (1,\infty)\setminus\mathbb{N}$, and any irrational $\alpha$.

\section{Examples}\label{SectionWithExamples}
Theorems \ref{GeneralizationOfClassicalvdC} and \ref{WeakMixingUDvdC} allows us to construct examples of sL-sequences and wm-sequences. Our next result shows us that we can use measure preserving systems to construct examples of sequences that appear in Definition \ref{MixingAndRigidDistributionsDefinitions}.
\begin{lemma}
    Let $\mathcal{X} := ([0,1]^d,\mathscr{B},m^d,T)$ be a measure preserving system, $x \in [0,1]^d$, and $x_n = T^nx$. The following hold for a.e. $x \in X$:
    \begin{enumerate}[(i)]
        \item If $\mathcal{X}$ is weakly mixing, then $(x_n)_{n = 1}^\infty$ is a wm-sequence.
        \item If $\mathcal{X}$ has Lebesgue spectrum, then $(x_n)_{n = 1}^\infty$ is a sL-sequence.
        \item If $\mathcal{X}$ has discrete spectrum, then $(x_n)_{n = 1}^\infty$ is a c-sequence. 
        \item If $\mathcal{X}$ has discrete spectrum and $A \in \mathscr{B}$, then $\{n_k\ |\ x_{n_k} \in A\}$ is a compact sequence of natural numbers.
        \item If $\mathcal{X}$ has singular spectrum, then $(x_n)_{n = 1}^\infty$ is a ss-sequence. 
        \item If $\mathcal{X}$ has singular spectrum and $A \in \mathscr{B}$, then $\{n_k\ |\ x_{n_k} \in A\}$ is a spectrally singular sequence of natural numbers.
    \end{enumerate}
\end{lemma}

\begin{proof}
    To prove (i) and (ii), we first observe that the system $\mathcal{X}$ is ergodic. Using the pointwise ergodic theorem, we see that for all $h \in \mathbb{N}$ and $f \in C([0,1]^d)\cup\{\mathbbm{1}_A\ |\ A \in \mathscr{B}\}$ we have

    \begin{equation}
        \lim_{N\rightarrow\infty}\frac{1}{N}\sum_{n = 1}^Nf(T^{n+h}x)\overline{f(T^nx)} = \int_0^1T^hf\overline{f}dm^d,
    \end{equation}
    for a.e. $x \in X$, from which the desired result follows. 

    We now proceed to prove items (iii)-(vi). Let $\{\mathcal{X}_u := (X_u,\mathscr{B}_u,m_u,T_u)\ |\ u \in U\}$ denote the ergodic components of $\mathcal{X}$. We note that if $\mathcal{X}$ has discrete (singular) spectrum, then $m^d$-a.e. $x \in [0,1]^d$ is contained in an ergodic component $\mathcal{X}_u$ that also has discrete (singular) spectrum. We observe that for $m_u$-a.e. $x \in X_u$, we have

    \begin{equation}\label{EquationUsingErgodicDecomposition}
        \gamma_f(h) := \lim_{N\rightarrow\infty}\frac{1}{N}\sum_{n = 1}^Nf(T^{n+h}x)\overline{f(T^nx)} = \int_XT^hf\overline{f}dm_u
    \end{equation}
    for all $f \in C([0,1]^d)$. Since $(\gamma_f(h))_{h = 1}^\infty$ are the Fourier coefficients of a measure that is absolutely continuous to the maximal spectral type of $\mathcal{X}_u$, we have proven (iii) and (v). To show (iv) and (vi), we repeat the previous argument with $f = \mathbbm{1}_A$ in Equation \eqref{EquationUsingErgodicDecomposition}.
\end{proof}

For a discussion of the prevalence of measure preserving systems with singular spectrum as well as examples, we refer the reader to \cite[Remark 1.16]{LebesgueSpectrumvanderCorput}. For now, we only mention one particularly aesthetic family of examples. Recall that the Thue-Morse sequence $(w_n)_{n = 1}^\infty \in \{-1,1\}^\mathbb{N}$ is given by $w_n = (-1)^{e_n}$ where $e_n$ is the number of times the digit $1$ appears in the base $2$ expansion of $n$. Let $(n_k)_{k = 1}^\infty$ be an increasing enumeration of those $n$ for which $w_n = 1$. The sequence $(n_k)_{k = 1}^\infty$ is a spectrally singular sequence as a consequence of the work of Mahler \cite{MahlerThue-Morse} (see also \cite{KakutaniThue-Morse}). Consequently, if $(x_n)_{n = 1}^\infty \subseteq [0,1]^d$ is a sL-sequence, then $(x_{n_k})_{k = 1}^\infty$ is uniformly distributed. 

More generally, for $r \in \mathbb{N}_{\ge 2}$, we define $s_r(n)$ to be the sum of the digits of the base $r$ expansion of $n$. It was shown in \cite{GeneralizedTMHasSS} that if $c \in \mathbb{R}$ is such that $(r-1)c \notin \mathbb{Z}$, then $((e(cs_r(n)))_{n = 1}^\infty,(e(cs_r(n)))_{n = 1}^\infty,(q)_{q = 1}^\infty)$ is a permissible triple, any that $(e(cs_r(n)))_{n = 1}^\infty$ is a spectrally singular sequence. We observe that $w_n = e\left(\frac{1}{2}s_2(n)\right)$. A consequence of this is the following. If $r \in \mathbb{N}_{\ge 2}$ and $0 \le j < r$, and $(n_k(j,r))_{k = 1}^\infty$ is an increasing enumeration of the $n \in \mathbb{N}$ for which $s_r(n) \equiv j\pmod{r}$, then $(n_k(j,r))_{k = 1}^\infty$ is a spectrally singular sequence. Consequently, if $(x_n)_{n = 1}^\infty \subseteq [0,1]^d$ is a sL-sequence, then $(x_{n_k(j,r)})_{k = 1}^\infty$ is uniformly distributed.

Our next example justifies the usage of permissible pairs in Theorem \ref{4CharacterizationsOfwmSequences}, and shows that Theorems \ref{GeneralizationOfClassicalvdC} and \ref{WeakMixingUDvdC} are not characterizations of sL-sequences and wm-sequences respectively.

\begin{lemma}\label{ExampleNecessitatingPermissibleTriples}
    There exists a sL-sequence $(x_n)_{n = 1}^\infty \subseteq [0,1]$ for which

    \begin{equation}\label{BadDiscrepancyEquation}
        \overline{D}((x_{n+h}-x_n)_{n = 1}^\infty) = \frac{1}{2}
    \end{equation}
    for all $h \in \mathbb{N}$.
\end{lemma}

\begin{proof}
    Let $(N_q)_{q = 0}^\infty$ be given by $N_0 = 0$ and $N_q = (q!)^2+N_{q-1}$. Let $f:\mathbb{N}\rightarrow\mathbb{N}^2$ the bijection corresponding to $\{(1,1),(1,2),(2,1),(3,1),(2,2),(1,3),(1,4),\cdots\}$, let $f_1(q)$ denote the first coordinate of $f(q)$, and observe that the first coordinate of $f(q)$ always divides $q!$. For all $q \in \mathbb{N}$, let us divide $(N_{q-1},N_q]$ into $M_q$ consecutive intervals of the form $\{(a_{i,q},a_{i,q}+2f_1(q)]\}$. We now define

    \begin{equation}
        x_n := \begin{cases}
                    n^2\sqrt{2}&\text{ if }n \in (a_{i,q},a_{i,q}+f_1(q)]\\
                    (n-f_1(q))^2\sqrt{2}&\text{ if }n \in (a_{i,q}+f_1(q),a_{i,q}+2f_1(q)].
               \end{cases}
    \end{equation}
    We will now show that Equation \eqref{BadDiscrepancyEquation} holds. Now let us fix some $h \in \mathbb{N}$ and let $(q_{h,k})_{k = 1}^\infty$ be an increasing enumeration of $f_1^{-1}(h)$. We see that for $n \in (N_{q_{h,k}-1},N_{q_{h,k}}-h]$

    \begin{equation}
        x_{n+h}-x_n = \begin{cases}
                    0&\text{ if }n \in (a_{i,q_{h,k}},a_{i,q_{h,k}}+h]\\
                    4nh\sqrt{2}&\text{ if }n \in (a_{i,q_{h,k}}+h,a_{i,q_{h,k}}+2h].
               \end{cases}
    \end{equation} 
    Now let

    \begin{equation}
        I_{1,k} = \bigcup_{i = 1}^{M_{q_{h,k}}}(a_{i,q_{h,k}},a_{i,q_{h,k}}+h]\text{ and }I_{2,k} = \bigcup_{i = 1}^{M_{q_{h,k}}}(a_{i,q_{h,k}}+h,a_{i,q_{h,k}}+2h].
    \end{equation}
    Firstly, we see that $D_{|I_{1,k}|}((x_{n+h}-x_n)_{n \in I_{1,k}}) = 1$. Secondly, we see that

    \begin{equation}
        (x_{n+h}-x_n)_{n \in I_{2,k}} = \bigcup_{j = h}^{2h}((N_{q_{h,k}-1}h+4hj+2hn')\sqrt{2})_{n' = 1}^{M_{q_{h,k}}}.
    \end{equation}
    Since $\sqrt{2}$ is an algebraic irrational, we may use Theorem \ref{DiscrepancyOfRotationByAnAlgebraicIrrational} to pick a constant $C = C(\sqrt{2},\frac{1}{2})$ such that for any $M < M_{q_{h,k}}$ we have

    \begin{equation}
        D_{M}((N_{q_{h,k}-1}h+4hj+2hn')\sqrt{2})_{n' = 1}^{M} \le CM^{-\frac{1}{2}}
    \end{equation}
    After recalling that the discrepancy function is subadditive, we see that

    \begin{alignat*}{2}
        &D_{|I_{2,k}|}((x_{n+h}-x_n)_{n \in I_{2,k}}) \le hCM_{q_{h,k}}^{-\frac{1}{2}}\text{, hence}\\
        &D_{N_{q_{h,k}}-N_{q_{h,k}-1}}\left((x_{n+h}-x_n)_{n = N_{q_{h,k}-1}+1}^{N_{q_{h,k}}}\right) \in \left(\frac{1}{2}-hCM_{q_{h,k}}^{-\frac{1}{2}},\frac{1}{2}+hCM_{q_{h,k}}^{-\frac{1}{2}}\right).
    \end{alignat*}
    Since $N_q >> N_{q-1}$, we see that

    \begin{equation*}
        \lim_{k\rightarrow\infty}D_{N_{q_{h,k}}}(x_{n+h}-x_n)_{n = 1}^{N_{q_{h,k}}} = \lim_{k\rightarrow\infty}D_{N_{q_{h,k}}-N_{q_{h,k}-1}}(x_{n+h}-x_n)_{n = N_{q_{h,k}-1}+1}^{N_{q_{h,k}}} = \frac{1}{2}.
    \end{equation*}
    The calculations thusfar show that $\overline{D}((x_{n+h}-x_n)_{n = 1}^\infty) \ge \frac{1}{2}$ for all $h \in \mathbb{N}$. For the reverse inequality, it suffices to show that
    
    \begin{equation}
        \lim_{\underset{q \notin f_1^{-1}(h)}{q\rightarrow\infty}}D_{N_q}((x_{n+h}-x_n)_{n = 1}^\infty) = 0.
    \end{equation}
    To this end, let $h_1 \neq h_2 \in \mathbb{N}$ be arbitrary, and observe that for\\ $n \in (N_{q_{h_1,k}-1},N_{q_{h_1,k}}-h_2]$ we have

    \begin{equation}
        x_{n+h_2}-x_n = ((n+f(n))^2-n^2)\sqrt{2} = 2f(n)n\sqrt{2}+f(n)^2\sqrt{2},
    \end{equation}
    where $f:\mathbb{N}\rightarrow\mathbb{Z}\setminus\{0\}$ is some periodic function with period $p$ dividing $2h_1$. Using the same bound as before for the discrepancy of $(n\sqrt{2})_{n = 1}^\infty$, we see that for $M \le M_{q_{h_1,k}}-h_2$ and $1 \le j \le 2h_1$ we have

    \begin{alignat*}{2}
        & D_M((x_{N_{q_{h_1,k}-1}+2h_1n'+j+h_2}-x_{N_{q_{h_1,k}-1}+2h_1n'+j})_{n' = 1}^M \le C(h_1,h_2,j)M^{-\frac{1}{2}}\text{, hence}\\
        &D_{2h_1M}((x_{N_{q_{h_1,k}-1}+n+h_2}-x_{N_{q_{h_1,k}-1}+n})_{n = 1}^{2h_1M}) \le C(h_1,h_2)M^{-\frac{1}{2}},
    \end{alignat*}
    from which the desired result follows.

    It remains to show that $(x_n)_{n = 1}^\infty$ is a sL-sequence. To this end, let $\ell \in \mathbb{N}$ be arbitrary and let $((e(\ell x_n))_{n = 1}^\infty,(e(\ell x_n))_{n = 1}^\infty,(N_m)_{m = 1}^\infty)$ be a permissible triple. It suffices to show that

    \begin{equation}\label{Provingx_nIssLEquation1}
        \lim_{m\rightarrow\infty}\frac{1}{N_q}\sum_{n = 1}^{N_m}e(\ell(x_{n+h}-x_n)) = 0,
    \end{equation}
    for all but at most $3$ values of $h \in \mathbb{N}$. If Equation \eqref{Provingx_nIssLEquation1} held for all $h \in \mathbb{N}$, then we would be done, so let $h_0 \in \mathbb{N}$ be such that

    \begin{equation}\label{Provingx_nIssLEquation2}
        \lim_{m\rightarrow\infty}\frac{1}{N_q}\sum_{n = 1}^{N_m}e(\ell(x_{n+h_0}-x_n)) = c \neq 0.
    \end{equation}
    Theorem \ref{SpecialCaseOfKHInequality} tells us that

    \begin{equation}
        \frac{1}{2}\cdot\frac{|c|}{4\pi \ell} \le D_{N_m}\left((x_{n+h_0}-x_n)_{n = 1}^{N_m}\right),
    \end{equation}
    for all $m \ge m_0$. Since $N_m >> N_{m-1}$, we must have that $N_{q_{h_0,k_m}} \le N_m < N_{q_{h_0,k_m}+1}$ for all $m \ge m_0$. Since $f_1(q_{h_0,k_m}+1) \in \{h_0-1,h_0,h_0+1\}$, we see that Equation \eqref{Provingx_nIssLEquation1} holds for all $h \notin \{h_0-1,h_0,h_0+1\}$.
\end{proof}

The sequence constructed in the previous Lemma is an example of an sL-sequence that is not an o-sequence, as any o-sequence $(x_n)_{n = 1}^\infty$ will have $(x_{n+h}-x_n)_{n = 1}^\infty$ be uniformly distributed for all $h \in \mathbb{N}$. Our next result gives an another example of a sL-sequence that is not an o-sequence. 

\begin{lemma}[{\cite[Theorem 2.4.28]{SohailsPhDThesis}}]
    The sequence $(x_n)_{n = 1}^\infty \subseteq [0,1]$ given by

    \begin{equation}
        x_n = \begin{cases}
                n^2\alpha&\text{ if }n\text{ is even}\\
                2(n-1)^2\alpha\text{ if }n\text{ is odd}
        \end{cases}
    \end{equation}
    is a sL-sequence but not an o-sequence. In particular, $(x_{n+h}-x_n)_{n = 1}^\infty$ is uniformly distributed for all $h \in \mathbb{N}$ and $(x_n,x_{n+h})_{n = 1}^\infty$ is uniformly distributed if and only if $h > 1$. 
\end{lemma}

\begin{acknowledgement}
I would like to thank Mariusz Lema\'nczyk for helpful discussions regarding spectral theory. I would also like to thank Teturo Kamae for reading an earlier draft of this paper and making some useful remarks. I acknowledge being supported by grant
2019/34/E/ST1/00082 for the project “Set theoretic methods in dynamics and number theory,” NCN (The
National Science Centre of Poland).
\end{acknowledgement}

\section{REFERENCES} 

\end{document}